\documentclass[11pt]{amsart}
\usepackage{mathrsfs}
\usepackage{amssymb,latexsym}
\usepackage{pstcol,pstricks,color}

 \numberwithin{equation}{section}

\newtheorem{theorem}{Theorem}[section]

\newtheorem{corollary}[theorem]{Corollary}

\newtheorem{remark}[theorem]{Remark}

\newtheorem{lemma}[theorem]{Lemma}

\def\qed{\hfill $\Box$}
\def\pf{\noindent {\it Proof.} }

\title{Congruences on the Bell polynomials and the derangement polynomials }

\begin{document}
\maketitle
\begin{center}
Yidong Sun$^\dag$,  Xiaojuan Wu and Jujuan Zhuang$^\ddag$


Department of Mathematics, Dalian Maritime University, 116026 Dalian, P.R. China\\[5pt]

{\it  Emails: $^\dag$sydmath@yahoo.com.cn, $^\ddag$jjzhuang1979@yahoo.com.cn}

\end{center}\vskip0.2cm

\subsection*{Abstract}
In this note, by the umbra calculus method, the Sun and Zagier's congruences
involving the Bell numbers and the derangement numbers are generalized to the polynomial cases.
Some special congruences are also provided.

\medskip

{\bf Keywords}:  Bell polynomials; Derangement polynomials; Stirling numbers; Congruences.

\noindent {\sc 2000 Mathematics Subject Classification}: Primary 11B75; Secondary 05A15, 05A18,
11A07.

\section{Introduction}

It is well known that the first and second kind Stirling numbers $s(m,j)$ and $S(m,j)$ \cite{Stanley} are defined respectively by
\begin{eqnarray}
x(x-1)\cdots (x-m+1)                     &=& \sum_{j=0}^{m}s(m,j)x^{j}, \label{eqn 1.1} \\
\sum_{j=0}^{m}S(m,j)x(x-1)\cdots (x-j+1) &=& x^m.  \label{eqn 1.2}
\end{eqnarray}

The Bell polynomials $\{\mathcal{B}_n(x)\}_{n\geq 0}$ are defined by
\begin{eqnarray*}
\mathcal{B}_m(x)=\sum_{j=0}^{m}S(m,j)x^j.
\end{eqnarray*}
It is clear that $\mathcal{B}_m(1)$ is the $m$-th Bell number, denoted by $B_m$, counting the number of partitions
of $[m]=\{1,2, \dots, m\}$ (with $B_{0} = 1$). The Bell polynomials $\mathcal{B}_m(x)$ satisfy the recurrence
\begin{eqnarray}\label{eqn 1.3}
\mathcal{B}_{m+1}(x)=x\sum_{j=0}^{m}\binom{m}{j}\mathcal{B}_{j}(x).
\end{eqnarray}

The derangement polynomials $\{\mathcal{D}_{m}(x)\}_{m\geq 0}$ are defined by
\begin{eqnarray*}
\mathcal{D}_{m}(x)=\sum_{j=0}^{m}\binom{m}{j}j!(x-1)^{m-j}.
\end{eqnarray*}
Clearly, $\mathcal{D}_m(1)=m!$ and $\mathcal{D}_m(0)$ is the $m$-th derangement number, denoted by $D_m$,
counting the number of fixed-point-free permutations on $[m]$ (with $D_{0}= 1$).
The derangement polynomials $\mathcal{D}_{m}(x)$, also called $x$-factorials of $m$, have been considerably
investigated by Eriksen, Freij and W$\ddot{a}$stlund \cite{Eriksen}, Sun and Zhuang \cite{SunZhuang}.
They obey the recursive relation
\begin{eqnarray}\label{eqn 1.4}
\mathcal{D}_{m}(x)=m\mathcal{D}_{m-1}(x)+(x-1)^{m}.
\end{eqnarray}

Recently, Sun \cite{Sun} discovered experimentally that for a fixed positive integer $m$
the sum $\sum_{k=0}^{p-1}B_{k}/(-m)^{k}$ modulo a prime $p$ not dividing $m$ is independent
of the prime $p$, a typical case being
\begin{eqnarray*}
\sum_{k=0}^{p-1}\frac{B_{k}}{(-8)^{k}} \equiv  -1853 \hskip.3cm  (mod\ p)\hskip.3cm  \mbox{ for\ all\ primes\ $p\neq 2$ }.
\end{eqnarray*}

Later, Sun and Zagier \cite{SunZag} confirmed this conjecture and proved the nice result.
\begin{theorem}
For any integer $m\geq 1$ and any prime $p\nmid m$, there hold
\begin{eqnarray*}
(-x)^{m}\sum_{k=1}^{p-1}\frac{\mathcal{B}_{k}(x)}{(-m)^{k}}
 \equiv   (-x)^p\sum_{k=0}^{m-1}\frac{(m-1)!}{k!}(-x)^{k}  \hskip.5cm (\mbox{mod}\ p).
\end{eqnarray*}
Particularly, the case $x=1$ generates
\begin{eqnarray}\label{eqn 1.5}
\sum_{k=1}^{p-1}\frac{B_{k}}{(-m)^{k}}  \equiv (-1)^{m-1}D_{m-1}  \hskip.5cm (\mbox{mod}\ p).
\end{eqnarray}
\end{theorem}

Here for two polynomials $P(x),Q(x) \in \mathbb{Z}_p[x]$, by $P(x)\equiv Q(x)\ (mod\ p)$ we mean that
the corresponding coefficients of $P(x)$ and $Q(x)$ are congruent modulo $p$.

In this note, we establish a more general result of Sun and Zagier's congruence.

\begin{theorem}
For any integers $n\geq 0, m\geq 1$ and any prime $p\nmid m$, there hold
\begin{eqnarray}\label{eqn 1.6}
x^{m}\sum_{k=1}^{p-1}\frac{\mathcal{B}_{n+k}(x)}{(-m)^{k}}
 \equiv   x^p\sum_{k=0}^{n}S(n,k)(-1)^{m+k-1}\mathcal{D}_{m+k-1}(1-x)  \hskip.5cm (\mbox{mod}\ p),
\end{eqnarray}
or equivalently
\begin{eqnarray}\label{eqn 1.7}
x^{m}\sum_{j=0}^{n}s(n,j)\sum_{k=1}^{p-1}\frac{\mathcal{B}_{j+k}(x)}{(-m)^{k}}
 \equiv    (-1)^{m+n-1}x^{p}\mathcal{D}_{m+n-1}(1-x)  \hskip.5cm (\mbox{mod}\ p).
\end{eqnarray}
\end{theorem}

In particular, the case $x=1$ leads to
\begin{corollary}
For any integers $n\geq 0, m\geq 1$ and any prime $p\nmid m$, there hold
\begin{eqnarray}\label{eqn 1.8}
\sum_{k=1}^{p-1}\frac{B_{n+k}}{(-m)^{k}}  \equiv  \sum_{k=0}^{n}S(n,k)(-1)^{m+k-1}D_{m+k-1}  \hskip.5cm (\mbox{mod}\ p),
\end{eqnarray}
or equivalently
\begin{eqnarray}\label{eqn 1.9}
\sum_{j=0}^{n}s(n,j)\sum_{k=1}^{p-1}\frac{B_{j+k}}{(-m)^{k}}  \equiv  (-1)^{m+n-1}D_{m+n-1}  \hskip.5cm (\mbox{mod}\ p).
\end{eqnarray}
\end{corollary}

\section{Proof of Theorem 1.2}

Define the generalized Bell umbra $\mathbf{B_x}$, given by $\mathbf{B_x}^{m}=\mathcal{B}_m(x)$.
(See \cite{Roman, RomRota} for more information on the umbra calculus.)
Then (\ref{eqn 1.3}) can be rewritten as $\mathbf{B_x}^{m+1}=x(\mathbf{B_x}+1)^{m}$. By linearity, for any polynomial $f(x)$ we have
\begin{eqnarray*}
\mathbf{B_x}f(\mathbf{B_x})=xf(\mathbf{B_x}+1),
\end{eqnarray*}
which, by induction on integer $m\geq 0$, yields
\begin{eqnarray}\label{eqn 2.1}
\mathbf{B_x}(\mathbf{B_x}-1)\cdots(\mathbf{B_x}-m+1)f(\mathbf{B_x})=x^{m}f(\mathbf{B_x}+m).
\end{eqnarray}

\begin{lemma} For any integers $m, n\geq 0$, there hold
\begin{eqnarray}\label{eqn 2.2a}
\mathbf{B_x}(\mathbf{B_x}-1)\cdots(\mathbf{B_x}-m+1)\mathbf{B_x}^{n}=x^{m}(\mathbf{B_x}+m)^{n}.
\end{eqnarray}
or equivalently
\begin{eqnarray}\label{eqn 2.2b}
\sum_{j=0}^{m}s(m,j)\mathcal{B}_{j+n}(x)= x^{m}\sum_{j=0}^{n}\binom{n}{j}\mathcal{B}_{j}(x)m^{n-j},
\end{eqnarray}
\end{lemma}
\pf The case $f(x)=x^{n}$ in (\ref{eqn 2.1}) produces (\ref{eqn 2.2a}).
By setting $x=\mathbf{B_x}$ in (\ref{eqn 1.1}), then (\ref{eqn 2.2a}) is just the umbral representation of (\ref{eqn 2.2b}). \qed

\begin{lemma} For any integer $m\geq 1$, there hold
\begin{eqnarray}\label{eqn 2.3a}
(\mathbf{B_x}-1)(\mathbf{B_x}-2)\cdots(\mathbf{B_x}-m+1) = (-1)^{m-1}\mathcal{D}_{m-1}(1-x),
\end{eqnarray}
or equivalently
\begin{eqnarray}\label{eqn 2.3b}
\sum_{j=0}^{m}s(m,j)\mathcal{B}_{j-1}(x) = (-1)^{m-1}\mathcal{D}_{m-1}(1-x).
\end{eqnarray}
\end{lemma}
\pf Let $\mathcal{A}_{m}(x)$ denote the expression on the left hand side of (\ref{eqn 2.3a}), by the case $n=0$ in (\ref{eqn 2.2a}), we have
\begin{eqnarray*}
\mathcal{A}_{m+1}(x)
&=& (\mathbf{B_x}-1)\cdots(\mathbf{B_x}-m+1)(\mathbf{B_x}-m) \\
&=& \mathbf{B_x}(\mathbf{B_x}-1)\cdots(\mathbf{B_x}-m+1)-m(\mathbf{B_x}-1)\cdots(\mathbf{B_x}-m+1) \\
&=& x^{m}-m\mathcal{A}_{m}(x).
\end{eqnarray*}
By (\ref{eqn 1.4}), it is routine to check that $(-1)^{m}\mathcal{D}_{m}(1-x)$ also obey the same recurrence as $\mathcal{A}_{m+1}(x)$ and
$\mathcal{A}_1(x)=1=\mathcal{D}_{0}(1-x)$. Hence $\mathcal{A}_{m+1}(x)=(-1)^{m}\mathcal{D}_{m}(1-x)$, which proves (\ref{eqn 2.3a}).

For (\ref{eqn 2.3b}), by setting $x=\mathbf{B_x}$ in (\ref{eqn 1.1})
after dividing an $x$ on the two sides of (\ref{eqn 1.1}), we can represent (\ref{eqn 2.3b}) umbrally as (\ref{eqn 2.3a})
and vice versa.  \qed

\begin{remark}
The case $x=1$ in (\ref{eqn 2.3b}) produces
\begin{eqnarray*}
\sum_{j=0}^{m}s(m,j)B_{j-1} = (-1)^{m-1}D_{m-1}.
\end{eqnarray*}
It is curious that such a simple and interesting identity did not
appear in the literature.
\end{remark}

\begin{remark}
By the orthogonal relationship between the two types of Stirling
numbers,
\begin{eqnarray}\label{eqn 2.4}
\sum_{j=k}^{m}s(m,j)S(j,k) = \delta_{m,k},
\end{eqnarray}
where $\delta_{m,k}$ is the Kronecker symbol defined by $\delta_{m,k}=1$ if $m=k$ and $\delta_{m,k}=0$ otherwise,
one can obtain another equivalent form of (\ref{eqn 2.2b}) and (\ref{eqn 2.3b})
\begin{eqnarray}
\mathcal{B}_{m+n}(x) &=& \sum_{k=0}^{m}S(m,k)x^{k}\sum_{j=0}^{n}\binom{n}{j}\mathcal{B}_{j}(x)k^{n-j}, \label{eqn 2.5} \\
\mathcal{B}_{m-1}(x) &=& \sum_{k=0}^{m}S(m,k)(-1)^{k-1}\mathcal{D}_{k-1}(1-x). \nonumber
\end{eqnarray}
It should be noticed that (\ref{eqn 2.5}) has been obtained by Spivey \cite{Spivey} in the case $x=1$, Gould
and Quaintance \cite{GouldQuain}, Belbachir and Mihoubi
\cite{BelMih} using different methods. By $S(p,1)=S(p,p)=1$ and $p|S(p,k)$ for a prime $p$ and $1<k<p$,
we have immediately the following congruence relations
\begin{eqnarray}
\mathcal{B}_{p+n}(x) &\equiv & x^{p}\mathcal{B}_{n}(x)+\mathcal{B}_{n+1}(x)  \hskip.5cm (\mbox{mod}\ p), \label{eqn 2.6} \\
\mathcal{B}_{p-1}(x) &\equiv & 1+\mathcal{D}_{p-1}(1-x) \hskip.5cm (\mbox{mod}\ p). \label{eqn 2.7}
\end{eqnarray}
Note that (\ref{eqn 2.6}) has been obtained by Gertsch and Robert \cite{GertRob},
and the case $x=1$ in (\ref{eqn 2.6}) reduces to the well-known Touchard's congruence
$B_{p+n} \equiv B_{n}+ B_{n+1}\ (\mbox{mod}\ p)$ \cite{Toucha}.
The case $x=1$ in (\ref{eqn 2.7}) yields a new congruence
\begin{eqnarray}\label{eqn 2.8}
B_{p-1} \equiv 1+D_{p-1} \hskip.5cm (\mbox{mod}\ p).
\end{eqnarray}
Later, Sun \cite{Sun2} informed us that they also independently obtained (\ref{eqn 2.8}) as a corollary of (\ref{eqn 1.5}).

\end{remark}\vskip.2cm

{\bf Proof of Theorem 1.2.} It suffices to prove (\ref{eqn 1.7}), for (\ref{eqn 1.6}) can be obtained
from (\ref{eqn 1.7}) by using the orthogonality in (\ref{eqn 2.4}). Setting $x=\mathbf{B_x}$ in the Lagrange congruence
\begin{eqnarray*}
x(x-1)\cdots (x-p+1)  \equiv x^{p}-x \hskip.5cm (\mbox{mod}\ p),
\end{eqnarray*}
by (\ref{eqn 2.2a}) in the case $n=0$, we have
\begin{eqnarray*}
\mathbf{B_x}^{p}-\mathbf{B_x}  \equiv x^{p} \hskip.5cm (\mbox{mod}\ p).
\end{eqnarray*}
Using the congruence $\binom{p-1}{k}\equiv (-1)^{p-k-1} (\mbox{mod}\ p)$ and
the Fermat's congruence $m^{p-1}\equiv 1\ (\mbox{mod}\ p)$, where $m$ is any integer not divided by the prime $p$, we get
\begin{eqnarray*}
\lefteqn{ x^{m}\sum_{j=0}^{n}s(n,j)\sum_{k=1}^{p-1}\frac{\mathcal{B}_{j+k}(x)}{(-m)^{k}} } \\
&=& x^{m}\sum_{k=1}^{p-1}\frac{1}{(-m)^{k}}\sum_{j=0}^{n}s(n,j)\mathcal{B}_{j+k}(x)        \\
&\equiv& x^{m}\sum_{k=1}^{p-1}\binom{p-1}{k}m^{p-k-1}x^{n}(\mathbf{B_x}+n)^{k}   \hskip.5cm (\mbox{mod}\ p)   \\
&=& x^{m+n}((\mathbf{B_x}+n+m)^{p-1}-m^{p-1})                                              \\
&=& \mathbf{B_x}(\mathbf{B_x}-1)\cdots(\mathbf{B_x}-(m+n)+1)(\mathbf{B_x}^{p-1}-m^{p-1})   \\
&\equiv& \mathbf{B_x}(\mathbf{B_x}-1)\cdots(\mathbf{B_x}-(m+n)+1)(\mathbf{B_x}^{p-1}-1) \hskip.5cm (\mbox{mod}\ p)  \\
&=& (\mathbf{B_x}-1)\cdots(\mathbf{B_x}-(m+n)+1)(\mathbf{B_x}^{p}-\mathbf{B_x})             \\
&\equiv& x^{p}(\mathbf{B_x}-1)\cdots(\mathbf{B_x}-(m+n)+1)   \hskip.5cm (\mbox{mod}\ p) \\
&=&  (-1)^{m+n-1}x^{p}\mathcal{D}_{m+n-1}(1-x),
\end{eqnarray*}
as desired. \qed

\section{ Special Consequences }

The cases $n=1$ and $n=2$ in (\ref{eqn 1.8}) produce
\begin{corollary} For any integer $m\geq 1$ and any prime $p\nmid m$, there hold
\begin{eqnarray}\label{eqn 3.1}
\sum_{k=1}^{p-1}\frac{B_{k+1}}{(-m)^{k}} &\equiv & (-1)^{m}D_{m}\hskip.5cm (\mbox{mod}\ p),          \\
\sum_{k=1}^{p-1}\frac{B_{k+2}}{(-m)^{k}}   &\equiv & (-1)^{m}(D_{m}-D_{m+1}) \hskip.5cm (\mbox{mod}\ p). \nonumber
\end{eqnarray}
\end{corollary}

\begin{corollary} For any integers $n, m\geq 0$ and any prime $p$, there hold
\begin{eqnarray}
D_{pn+m}                                        &\equiv &  (-1)^{n}D_{m} \hskip.5cm (\mbox{mod}\ p),   \label{eqn 3.2} \\
\sum_{k=1}^{p-1}(-1)^{k}B_{n+k}                 &\equiv &  V_{n} \hskip.5cm (\mbox{mod}\ p) ,        \label{eqn 3.3}  \\
\sum_{k=1}^{p-1}B_{n+k} -\sum_{k=1}^{n-1}B_{k}  &\equiv &  D_{p-1} \hskip.5cm (\mbox{mod}\ p).       \label{eqn 3.4}
\end{eqnarray}
where $V_n=\sum_{k=0}^{n}(-1)^{n-k}\binom{n}{k}B_k$ is the number of partitions of $[n]$ without singletons (i.e., one-element subsets) \cite{SunWu}.
\end{corollary}
\pf Setting $m:=pn+m$ in (\ref{eqn 3.1}), by $(-1)^{pn}\equiv (-1)^n\ (mod\ p)$, one can get (\ref{eqn 3.2}).
Setting $m=1$ in (\ref{eqn 1.8}), we have
\begin{eqnarray*}
\sum_{k=1}^{p-1}(-1)^{k}B_{n+k}
&\equiv & \sum_{k=0}^{n}S(n,k)(-1)^{k}D_{k}  \hskip.5cm (\mbox{mod}\ p) \\
&=& \sum_{k=0}^{n}S(n,k)(\mathbf{B}-1)(\mathbf{B}-2)\cdots(\mathbf{B}-k)  \\
&=& (\mathbf{B}-1)^{n}=\sum_{k=0}^{n}(-1)^{n-k}\binom{n}{k}\mathbf{B}^{k} \\
&=& \sum_{k=0}^{n}(-1)^{n-k}\binom{n}{k}B_k=V_n,
\end{eqnarray*}
where $\mathbf{B}:=\mathbf{B_x}|_{x=1}$ is the Bell umbra. Thus (\ref{eqn 3.3}) follows.

By the Touchard's congruence, we have
\begin{eqnarray*}
\lefteqn{ \sum_{k=1}^{p-1}B_{n+k} -\sum_{k=1}^{n-1}B_{k} }  \\
&\equiv& \sum_{k=1}^{p-1}B_{n+k}-\sum_{j=1}^{n-1}(B_{p+j}-B_{j+1} )  \hskip.5cm (\mbox{mod}\ p) \\
&=& \sum_{k=1}^{p-1}B_{n+k}-\sum_{j=1}^{n-1}\sum_{k=1}^{p-1}(B_{j+k+1}-B_{j+k} ) \\
&=& \sum_{k=1}^{p-1}B_{n+k}-\sum_{k=1}^{p-1}\sum_{j=1}^{n-1}(B_{j+k+1}-B_{j+k} ) \\
&=& \sum_{k=1}^{p-1}B_{n+k}-\sum_{k=1}^{p-1}(B_{n+k}-B_{k+1} ) \\
&=& \sum_{k=1}^{p-1}B_{k+1} \equiv   D_{p-1} \hskip.5cm (\mbox{mod}\ p),
\end{eqnarray*}
where the last step is obtained by setting $m:=p-1$ in (\ref{eqn 3.1}). \qed

\section*{Acknowledgements} { The authors are grateful to Prof. Zhi-Wei Sun
for the helpful suggestions and comments. The work was
partially supported by The National Science Foundation of China
(Grant No. 10801020 and 11001035) and by the Fundamental Research Funds
for the Central Universities (Grant No. 2009QN070 and 2009QN071).}


\end{document}